\documentclass[12pt]{article}
\usepackage{color}
\definecolor{darkblue}{rgb}{0.00,0.25,0.50}
\usepackage[colorlinks,filecolor=blue,citecolor=darkblue]{hyperref}

\usepackage{amsfonts,amssymb,amsmath,amsthm}
\usepackage{url}
\usepackage{enumerate}
\usepackage[ukrainian, russian, english]{babel}
\usepackage[cp1251]{inputenc}
\begin{document}\selectlanguage{ukrainian}
\thispagestyle{empty}

\title{}

\begin{center}
\textbf{\Large Порядкові оцінки найкращих ортогональних тригонометричних наближень класів згорток періодичних функцій невеликої  гладкості}
\end{center}
\vskip0.5cm
\begin{center}
А.~С.~Сердюк${}^1$, Т.~А.~Степанюк${}^2$\\ \emph{\small
${}^1$Інститут математики НАН
України, Київ\\
${}^2$Східноєвропейський національний університет імені Лесі
Українки, Луцьк\\}
\end{center}
\vskip0.5cm


\begin{abstract}
Получено порядковые оценки для наилучших равномерных ортогональных тригонометрических приближений на классах
$2\pi$--периодических функций, таких, что их
$(\psi,\beta)$--производные принадлежат единичным шарам пространств  $L_{p}, \ 1\leq p<\infty$, в случае когда последовательность
$\psi$ такая, что произведение $\psi(n)n^{\frac{1}{p}}$, $1\leq p<\infty$,
может стремится к нулю медленее за любую степенную функцию и $\sum\limits_{k=1}^{\infty}\psi^{p'}(k)k^{p'-2}<\infty$ при $1<p<\infty$, $\frac{1}{p}+\frac{1}{p'}=1$ или $\sum\limits_{k=1}^{\infty}\psi(k)<\infty$
 при $p=1$.
Аналогичные оценки получены для приближений  в $L_{s}$--метриках, $1<
s\leq \infty$, для классов  $(\psi,\beta)$--дифференцируемых
функций, таких, что $\parallel f_{\beta}^{\psi}\parallel_{1}\leq1$.

\vskip 0.5cm

 We obtain order estimates for the best uniform orthogonal trigonometric
 approximations of  $2\pi$--periodic
 functions, whose $(\psi,\beta)$--derivatives
 belong to unit balls of spaces $L_{p}, \ 1\leq
p<\infty$, in case at consequences  $\psi(k)$ are that product $\psi(n)n^{\frac{1}{p}}$ can tend to zero slower than any power function and
$\sum\limits_{k=1}^{\infty}\psi^{p'}(k)k^{p'-2}<\infty$ when $1<p<\infty$, $\frac{1}{p}+\frac{1}{p'}=1$ and $\sum\limits_{k=1}^{\infty}\psi(k)<\infty$
when $p=1$.
 We also established the analogical
estimates in $L_{s}$--metric, $1< s\leq \infty$, for classes
 of the summable
$(\psi,\beta)$--differentiable functions, such that  $\parallel
f_{\beta}^{\psi}\parallel_{1}\leq1$.

\end{abstract}

\vskip 1.5cm


Позначимо через
$L_{p}$,
$1\leq p<\infty$, --- простір $2\pi$--періодичних сумовних в $p$--му
степені на $[0,2\pi)$  функцій $f:\mathbb{R}\rightarrow\mathbb{C}$ з нормою
$${\|f\|_{p}:=\Big(\int\limits_{0}^{2\pi}|f(t)|^{p}dt\Big)^{\frac{1}{p}}},$$ а через
$L_{\infty}$ --- простір
$2\pi$--періодичних вимірних і суттєво обмежених  функцій ${f:\mathbb{R}\rightarrow\mathbb{C}}$ з нормою
$$\|f\|_{\infty}:=\mathop{\rm{ess}\sup}\limits_{t}|f(t)|.$$

Розглянемо множини $2\pi$--періодичних дійснозначних функцій $L^{\psi}_{\beta}$, які означаються наступним чином.

Нехай $f:\mathbb{R}\rightarrow\mathbb{R}$ --- функція із $L_{1}$, ряд Фур'є  якої має вигляд
$$
\sum_{k=-\infty}^{\infty}\hat{f}(k)e^{ikx},
$$
де
\begin{equation}\label{Fourier}
\hat{f}(k)=\frac{1}{2\pi}\int\limits_{-\pi}^{\pi}f(t)e^{-ikt}dt.
\end{equation}

\noindent Нехай, далі, $\psi(k)$ --- довільна фіксована послідовність дійсних чисел і $\beta$
--- фіксоване  дійсне число. Тоді якщо ряд
$$
\sum_{k\in \mathbb{Z}/\{0\}}\frac{\hat{f}(k)}{\psi(|k|)}e^{i(kx+\frac{\beta\pi}{2} \mathrm{sign} k)}
$$
\noindent є рядом Фур'є деякої сумовної функції $\varphi$ із $L_{1}$, то цю
функцію  називають (див., наприклад, \cite[с.
132]{Stepanets1}) $(\psi,\beta)$-похідною функції $f$ і позначають через
$f_{\beta}^{\psi}$.
Множину функцій $f$, у яких існує $(\psi,\beta)$-похідна
позначають через $L_{\beta}^{\psi}$.

Розглянемо одиничну кулю $B_{p}$ в просторі дійснозначних функцій з $L_{p}$, тобто множину функцій $\varphi:\mathbb{R}\rightarrow\mathbb{R}$ таких, що
$\|\varphi\|_{p}\leq1, \ 1\leq p\leq\infty$.
Якщо $f\in L^{\psi}_{\beta}$ і водночас $f^{\psi}_{\beta}\in
B_{p}$, то будемо записувати, що функція $f\in L^{\psi}_{\beta,p}$.

Як показано в \cite[с.
136]{Stepanets1}, якщо послідовність $\psi(k)$
  монотонно прямує до нуля при  $k\rightarrow\infty$
  i $\sum\limits_{k=1}^{\infty}\frac{\psi(k)}{k}<\infty$, то елементи $f(x)$ множини
  $L^{\psi}_{\beta,p}$,  $\beta\in\mathbb{R}$ майже  при всіх $x\in\mathbb{R}$ можна зобразити у вигляді згортки
 \begin{equation}\label{psi,beta}
f(x)=\frac{a_{0}}{2}+\frac{1}{\pi}\int\limits_{-\pi}^{\pi}\Psi_{\beta}(x-t)\varphi(t)dt,
\ a_{0}\in\mathbb{R}, \  \varphi\in B_{p}, \ \varphi\perp 1,
\end{equation}
з сумовним ядром $\Psi_{\beta}$,  ряд Фур'є якого має вигляд
 $$
\frac{1}{2}\sum\limits_{\mathbb{Z}/\{0\}}\psi(|k|)e^{-i(kt+\frac{\beta\pi}{2}\mathrm{sign} k)}=\sum\limits_{k=1}^{\infty}\psi(k)\cos
\big(kt-\frac{\beta\pi}{2}\big).
$$
При цьому функція $\varphi$ майже скрізь  збігається з
$f^{\psi}_{\beta}$.

Будемо вважати, що послідовності $\psi (k),\ k\in \mathbb{N},$ які задають  класи
 $L^{\psi}_{\beta,p}$, є звуженнями на множину
натуральних чисел  деяких додатних, неперервних, опуклих
донизу функцій $\psi(t)$, заданих на $[1, \infty)$, що задовольняють умову $
\lim\limits_{t\rightarrow\infty}\psi(t)=0. $  Множину всіх таких
функцій $\psi$ позначатимемо через ${\mathfrak M}$.

Для класифікації функцій $\psi$ із $\mathfrak{M}$ за їх швидкістю спадання до нуля важливу роль відіграє характеристика
\begin{equation}\label{for301}
\alpha(\psi;t):=\frac{\psi(t)}{t|\psi'(t)|}, \ \ \psi'(t):=\psi'(t+0).
\end{equation}
З її допомогою з множини ${\mathfrak M}$ виділяють наступні підмножини (див., наприклад, \cite[с.
160--161]{Stepanets1}):
\begin{equation}\label{m0}
  \mathfrak{M}_{0}:=\left\{\psi\in \mathfrak{M}: \ \ \exists K>0 \ \ \ \  \forall t\geq1 \ \ \ \
0<K\leq \alpha(\psi;t) \right\},
\end{equation}
\begin{equation}\label{mc}
  \mathfrak{M}_{C}:=\left\{\psi\in \mathfrak{M}: \ \ \exists K_{1}, K_{2}>0 \ \ \ \ \forall t\geq1 \ \ \ \
K_{1}\leq \alpha(\psi;t)\leq K_{2}<\infty \right\}.
\end{equation}
В (\ref{m0}) i (\ref{mc}) величини $K, \ K_{1}, \ K_{2}$ можуть залежати від $\psi$.
Очевидно, що $\mathfrak{M}_{C}\subset\mathfrak{M}_{0}$.

Нехай $m\in\mathbb{N}$, $\gamma_{m}$ --- довільний набір із $m$ цілих чисел і
$$
S_{\gamma_{m}}(f;x)=\sum\limits_{k\in \gamma_{m}}\hat{f}(k)e^{ikx},
$$
де  $\hat{f}(k)$ --- коефіцієнти Фур'є функції $f$ вигляду (\ref{Fourier}).

Величину
\begin{equation}\label{n_term1}
  e^{\bot}_{m}(f)_{s}=\inf\limits_{\gamma_{m}}\|f(x)-S_{\gamma_{m}}(f;x) \|_{s},\ 1\leq s\leq\infty,
\end{equation}
називають найкращим ортогональним тригонометричним наближенням функції ${f\in L_{s}}$ в метриці простору $L_{s}$,
а величину
\begin{equation}\label{n_term}
  e^{\bot}_{m}(L_{\beta,p}^{\psi})_{s}=\sup\limits_{f\in L_{\beta,p}^{\psi}}e^{\bot}_{m}(f)_{s}, \ 1\leq p,s\leq\infty,
\end{equation}
--- найкращим ортогональним тригонометричним наближенням класу $L_{\beta,p}^{\psi}$ в метриці простору $L_{s}$.

Метою даної роботи  є знаходження точних порядкових оцінок величин $ e^{\bot}_{n}(L_{\beta,p}^{\psi})_{s},  \ \beta\in\mathbb{R}$, при $ 1\leq p<\infty$ і $s=\infty$, а також при $p=1$ i $1< s\leq\infty$.

У випадку коли $\psi(k)=k^{-r},\ r>0$, класи $L_{\beta,p}^{\psi}, \ 1\leq p\leq\infty, \ \beta\in\mathbb{R}$ є
відомими класами
 Вейля-Надя $W_{\beta,p}^r$.
Для цих класів порядкові оцінки величин $(\ref{n_term})$ при $1<p,s<\infty$,  відомі (див. \cite{Romanyuk2002}, \cite{Romanyuk2012}).
Точні порядки величин $e^{\bot}_{n}(W^{r}_{\beta,p})_{s}$, $\beta\in\mathbb{R}$,  встановлені
також при $1< p<\infty$, $s=\infty$ для усіх  $r>\frac{1}{p}$, при $p=1$, $1<s<\infty$ для всіх $r>\frac{1}{s'}$, та при $s=\infty$, $p=1$, $r>1$ і $\beta=0$ 
(див. \cite{Romanyuk2007},  \cite[с. 137, 140]{Romanyuk2012}).

У випадку, коли $\psi\in B\cap \Theta_{p}^{*}$, де  $B$ --- множина
незростаючих додатних функцій $\psi(t)$, $t\geq 1$, для кожної з
яких можна вказати додатну сталу $K$ таку, що $
\frac{\psi(t)}{\psi(2t)}\leq K, \ \  t\geq 1 $, а $\Theta_{q}^{*}$ --- множина
незростаючих додатних функцій $\psi(t)$, для яких існує $\varepsilon>0$ таке, що послідовність $\psi(k)k^{\frac{1}{q}+\varepsilon}$ не зростає, в \cite{Shkapa1} були знайдені точні порядкові оцінки величин $e^{\bot}_{n}(L_{\beta,p}^{\psi})_{\infty}, \ 1< p<\infty$, $\beta\in\mathbb{R}$.
Якщо ж $\psi\in B\cap \Theta_{s'}^{*}$ і $\frac{1}{\psi(t)}$ опукла, то в
 роботі \cite{Shkapa}
  встановлені точні порядкові оцінки величин $e^{\bot}_{n}(L_{\beta,1}^{\psi})_{s}, \ 1< s<\infty$, для довільних $\beta\in\mathbb{R}$.

Зазначимо, що при довільних $1<p,s<\infty$ і $\beta\in\mathbb{R}$ точні порядки величин $e^{\bot}_{n}(L_{\beta,p}^{\psi})_{s}$ також відомі  (див., наприклад, \cite{Fedorenko1999} i \cite{Fedorenko2000}).

В даній роботі знайдено двосторонні  оцінки для
 величин $ e^{\bot}_{n}(L_{\beta,p}^{\psi})_{\infty}$,  $ 1\leq p<\infty $, у випадку, коли  функція
 $g_{p}(t)=\psi(t)t^{\frac{1}{p}}$ належить до множини $\mathfrak{M}_{0}$
i крім цього $${\left\{\begin{array}{cc}
\sum\limits_{k=1}^{\infty}\psi^{s}(k)k^{s-2}<\infty,\ \ & 1<s<\infty, \\
\sum\limits_{k=1}^{\infty}\psi(k)<\infty,  \ \ \ \ \  & s=1. \
  \end{array} \right.}$$
В ній також знайдено двосторонні  оцінки для
 величин  $ e^{\bot}_{n}(L_{\beta,1}^{\psi})_{s}$  у випадку коли
 ${g_{s'}\in\mathfrak{M}_{0}}$
i $\sum\limits_{k=1}^{\infty}\psi^{s}(k)k^{s-2}<\infty$, $1<s<\infty$, $\frac{1}{s}+\frac{1}{s'}=1$.
 При цьому константи в отриманих оцінках будуть виражені через параметри  класів в явному вигляді.

Позначимо через ${\cal E}_{n}(L^{\psi}_{\beta,p})_{s}$ наближення сумами Фур'є класів $L^{\psi}_{\beta,p}$ в метриках просторів $L_{s}$,
тобто величини
  вигляду
  \begin{equation}\label{sum}
  {\cal E}_{n}(L^{\psi}_{\beta,p})_{s}=\sup\limits_{f\in
L^{\psi}_{\beta,p}}\|f(\cdot)-S_{n-1}(f;\cdot)\|_{s}, \ 1\leq p,s\leq\infty,
\end{equation}
де $S_{n-1}(f;\cdot)$ --- частинні суми Фур'є порядку $n-1$ функції $f$.

З означень величин (\ref{n_term}) i (\ref{sum})
є очевидною нерівність
\begin{equation}\label{ineq_comp}
 e^{\bot}_{2n-1}(L_{\beta,p}^{\psi})_{s}\leq {\cal E}_{n}(L^{\psi}_{\beta,p})_{s}, \ 1\leq p,s\leq\infty.
\end{equation}

Отже, величини $ {\cal E}_{n}(L^{\psi}_{\beta,p})_{s}$ природньо використовувати для оцінки зверху найкращих ортогональних тригонометричних наближень
вигляду (\ref{n_term}). Встановленню точних порядкових оцінок величин $ {\cal E}_{n}(L^{\psi}_{\beta,p})_{s}$ при $1\leq p<\infty$ i $s=\infty$ та $p=1$ i $1<s\leq\infty$ присвячено роботи \cite{Serdyuk_grabova}-\cite{S_S}.

Щоб сформулювати основні результати роботи введемо наступні позначення. Для кожного $1<s<\infty$ покладемо
\begin{equation}\label{for100}
\xi(s):=\max\Big\{4\Big(\frac{\pi}{s-1}\Big)^{\frac{1}{s}}, \ \ 14(8\pi)^{\frac{1}{s}} s\Big\},
\end{equation}
а для будь--якої функції $\psi\in\mathfrak{M}$ через $\underline{\alpha}_{n}(\psi)$ i $\overline{\alpha}_{n}(\psi)$, $n\in \mathbb{N}$,
будемо позначати величини
\begin{equation}\label{k}
\underline{\alpha}_{n}(\psi):=\inf\limits_{ t\geq n}\alpha(\psi;t),
\end{equation}
\begin{equation}\label{kk}
\overline{\alpha}_{n}(\psi):=\sup\limits_{ t\geq n}\alpha(\psi;t),
\end{equation}
де характеристика $\alpha(\psi;t)$ означається формулою (\ref{for301}).
В прийнятих позначеннях має місце наступне твердження.

{\bf Теорема 1.} {\it Нехай $1<p<\infty, \ \sum\limits_{k=1}^{\infty}\psi^{p'}(k)k^{p'-2}<\infty$, $\frac{1}{p}+\frac{1}{p'}=1$,
 а функція ${g_{p}(t)=\psi(t)t^{\frac{1}{p}}}$ така, що
 $$
 g_{p}\in\mathfrak{M}_{0}
 $$
  і
  $$
  \underline{ \alpha}_{1}(g_{p})=\inf\limits_{t\geq1}\alpha(g_{p};t)>p'.
  $$
   Тоді для довільних $n\in \mathbb{N}$ i $\beta\in\mathbb{R}$ мають місце співвідношення
 \begin{equation}\label{theorem_1}
K^{(1)}_{\psi,p}\Big(\sum\limits_{k=n}^{\infty}\psi^{p'}(k)k^{p'-2}\Big)^{\frac{1}{p'}}\leq e^{\bot}_{2n}(L_{\beta,p}^{\psi})_{\infty}\leq e^{\bot}_{2n-1}(L_{\beta,p}^{\psi})_{\infty}\leq
K^{(2)}_{\psi,p}\Big(\sum\limits_{k=n}^{\infty}\psi^{p'}(k)k^{p'-2}\Big)^{\frac{1}{p'}},
\end{equation}
в яких
\begin{equation}\label{k1}
  K_{\psi,p}^{(1)}=\frac{1}{3\xi(p)}\Big(\frac{\underline{\alpha}_{1}(g_{p})}{p'+\underline{\alpha}_{1}(g_{p})}\Big)^{\frac{1}{p}}
\Big(1-\frac{p'}{\underline{\alpha}_{1}(g_{p})}\Big),
\end{equation}
\begin{equation}\label{k2}
  K_{\psi,p}^{(2)}=\frac{1}{\pi}\xi(p')\Big(\frac{p'+\underline{\alpha}_{1}(g_{p})}{\underline{\alpha}_{1}(g_{p})}\Big)^{\frac{1}{p'}}.
\end{equation}
}

{\bf Доведення теореми 1.}
Згідно з теоремою 1 роботи  \cite{Serdyuk_Stepaniuk2015} при виконанні умов ${\psi(t)t^{\frac{1}{p}}\in\mathfrak{M}_{0}}$ i
$\sum\limits_{k=1}^{\infty}\psi^{p'}(k)k^{p'-2}<\infty$, $1<p<\infty$, $\frac{1}{p}+\frac{1}{p'}=1$,
   $n\in \mathbb{N}$, $\beta\in\mathbb{R}$ справедлива оцінка
\begin{equation}\label{th1}
  {\cal E}_{n}(L^{\psi}_{\beta,p})_{\infty}
\leq
K_{\psi,p}^{(2)}\Big(\sum\limits_{k=n}^{\infty}\psi^{p'}(k)k^{p'-2}\Big)^{\frac{1}{p'}},
\end{equation}
в якій величини  $K_{\psi,p}^{(2)}$ означені формулою (\ref{k2}).

Враховуючи нерівності (\ref{ineq_comp}) i (\ref{th1}), отримуємо
\begin{equation}\label{w1}
e^{\bot}_{2n}(L_{\beta,p}^{\psi})_{\infty}\leq e^{\bot}_{2n-1}(L_{\beta,p}^{\psi})_{\infty}\leq
K^{(2)}_{\psi,p}\Big(\sum\limits_{k=n}^{\infty}\psi^{p'}(k)k^{p'-2}\Big)^{\frac{1}{p'}}.
\end{equation}
Знайдемо  оцінку знизу величини $e^{\bot}_{2n}(L_{\beta,p}^{\psi})_{\infty}$.
Розглянемо функцію
\begin{equation}\label{eq4}
  f^{*}_{p}(t)= f^{*}_{p}(\psi;n;t):=
\frac{\lambda}{\Big(\sum\limits_{k=n}^{\infty}\psi^{p'}(k)k^{p'-2}\Big)^{\frac{1}{p}}}\sum\limits_{k=n}^{\infty}
\psi^{p'}(k)k^{p'-2}\cos
kt,
\end{equation}
де
\begin{equation}\label{a}
 \lambda=\lambda(\psi;p;n):=\frac{1}{\xi(p)}\Big(\frac{\underline{\alpha}_{n}(g_{p})}{p'+\underline{\alpha}_{n}(g_{p})}\Big)^{\frac{1}{p}},  \ 1<p<\infty,
\ \ \frac{1}{p}+\frac{1}{p'}=1.
\end{equation}
 В \cite{Serdyuk_Stepaniuk2015} було показано, що при виконанні умови
 $g_{p}\in\mathfrak{M}_{0}$ функція $ f^{*}_{p}$ належить до $L_{\beta,p}^{\psi}, \ {1<p<\infty}$.
 Покажемо, що
 \begin{equation}\label{n10}
e^{\bot}_{2n}(f^{*}_{p})_{\infty}  \geq K^{(1)}_{\psi,p}\Big(\sum\limits_{k=n}^{\infty}\psi^{p'}(k)k^{p'-2}\Big)^{\frac{1}{p'}}.
\end{equation}

Нехай
\begin{equation}\label{Phi_s}
\Phi_{s}(x):=\int\limits_{x}^{\infty}\psi^{s}(t)t^{s-2}dt
\end{equation}
  і
\begin{equation}\label{for3331}
A_{s}(l;n)=A_{s}(\psi;l;n):=\big[\Phi_{s}^{-1}\big(\frac{1}{2l}\Phi_{s}(n)\big)\big]+2n, \ l\in\mathbb{N},
\end{equation}
 де $[\alpha]$ --- ціла частина дійсного числа $\alpha$,  $\Phi_{s}^{-1}$ --- функція обернена до $\Phi_{s}$.

Розглянемо величину
\begin{equation}\label{eq3}
  I_{1}:=\inf\limits_{\gamma_{2n}}\bigg|\int\limits_{-\pi}^{\pi}(f^{*}_{p}(t)-S_{\gamma_{2n}}(f^{*}_{p};t))V_{A_{p'}(l;n)}(t)dt\bigg|,
\end{equation}
де  $V_{A_{p'}(l;n)}$ --- ядра Валле Пуссена $V_{m}$ (див., наприклад,
 \cite[с. 31]{Stepanets1})
\begin{equation}\label{val_pus2}
  V_{m}(t)=\frac{1}{2}+\sum\limits_{k=1}^{m}\cos kt+2\sum\limits_{k=m+1}^{2m-1}\Big(1-\frac{k}{2m}\Big)\cos
kt, \ m\in
\mathbb{N},
\end{equation}
при $m=A_{p'}(l;n)$.

 В силу твердження Д.1.1 з 
 \cite[с. 391]{Korn} 
\begin{equation}\label{for200}
 I_{1}\leq\inf\limits_{\gamma_{2n}}\|f^{*}_{p}(t)-S_{\gamma_{2n}}(f^{*}_{p};t)\|_{\infty}\|V_{A_{p'}(l;n)}\|_{1}=
 e^{\bot}_{2n}(f^{*}_{p})_{\infty}\|V_{A_{p'}(l;n)}\|_{1}.
\end{equation}
Оскільки (див., наприклад, \cite[с.
247]{Stepaniuk2014})
\begin{equation}\label{eq61}
 \|V_{m}\|_{1}\leq3\pi, \ \ m\in\mathbb{N},
\end{equation}
то з (\ref{for200}) і (\ref{eq61}) можемо записати оцінку
\begin{equation}\label{eq9}
 e^{\bot}_{2n}(f^{*}_{p})_{\infty} \geq \frac{1}{3\pi}I_{1}.
\end{equation}
Ядра $V_{m}$ вигляду (\ref{val_pus2}) можна зобразити у вигляді
$$
V_{m}(t)=\frac{1}{2}\Big(1+\sum\limits_{1\leq k\leq m}e^{ikt}+\sum\limits_{-m\leq k\leq -1}e^{ikt}+
$$
\begin{equation}\label{val_pus_complex}
+2\sum\limits_{m+1\leq k\leq 2m-1}\Big(1-\frac{k}{2m}\Big)e^{ikt}
+2\sum\limits_{-2m+1\leq k\leq-m-1}\Big(1-\frac{|k|}{2m}\Big)e^{ikt}\Big).
\end{equation}

Крім того
 \begin{equation}\label{n1}
 f^{*}_{p}(t)-S_{\gamma_{2n}}(f^{*}_{p};t)=
\frac{\lambda}{2\Big(\sum\limits_{k=n}^{\infty}\psi^{p'}(k)k^{p'-2}\Big)^{\frac{1}{p}}}
{\mathop{\sum}\limits_{
 |k|\geq n,\atop k \notin\gamma_{2n} }}
\psi^{p'}(|k|)|k|^{p'-2}e^{ikt}.
\end{equation}
Оскільки
 \begin{equation}\label{int_riv}
  \int\limits_{-\pi}^{\pi}e^{ikt}e^{imt}dt=
{\left\{\begin{array}{cc}
0, \ & k+m\neq 0, \\
2\pi, &
k+m=0, \
  \end{array} \right.}  \ \ k,m\in\mathbb{Z},
\end{equation}
то з урахуванням (\ref{val_pus_complex}) маємо
$$
\int\limits_{-\pi}^{\pi}{\mathop{\sum}\limits_{
 |k|\geq n,\atop k \notin\gamma_{2n} }}
\psi^{p'}(|k|)|k|^{p'-2}e^{ikt}V_{A_{p'}(l;n)}(t)dt=
$$
$$
=\frac{1}{2}\int\limits_{-\pi}^{\pi}\bigg({\mathop{\sum}\limits_{
 k\geq n,\atop k \notin\gamma_{2n} }}
\psi^{p'}(k)k^{p'-2}e^{ikt}+{\mathop{\sum}\limits_{
 k\leq -n,\atop k \notin\gamma_{2n} }}
\psi^{p'}(|k|)|k|^{p'-2}e^{ikt}\bigg)\times
$$
$$
\times\Big(1+\sum\limits_{1\leq k\leq A_{p'}(l;n)}e^{ikt}+\sum\limits_{-A_{p'}(l;n)\leq k\leq-1}e^{ikt}+
$$
$$
+
2\sum\limits_{A_{p'}(l;n)+1\leq k\leq 2A_{p'}(l;n)-1}\Big(1-\frac{k}{2A_{p'}(l;n)}\Big)e^{ikt}+
$$
$$
+2\sum\limits_{-2A_{p'}(l;n)+1\leq k\leq-A_{p'}(l;n)-1}\Big(1-\frac{|k|}{2A_{p'}(l;n)}e^{ikt}\Big)\Big) dt =
$$
$$
=\pi\bigg(\mathop{\sum}\limits_{
n\leq k\leq A_{p'}(l;n),\atop k \notin\gamma_{2n} }
\psi^{p'}(k)k^{p'-2}+
\mathop{\sum}\limits_{
-A_{p'}(l;n)\leq k\leq-n,\atop k \notin\gamma_{2n} }
\psi^{p'}(|k|)|k|^{p'-2}+
$$
$$
+2
\mathop{\sum}\limits_{
A_{p'}(l;n)+1\leq k\leq 2A_{p'}(l;n)-1,\atop k \notin\gamma_{2n} }
\Big(1-\frac{k}{2A_{p'}(l;n)}\Big)\psi^{p'}(k)k^{p'-2}
+
$$
$$
+2\mathop{\sum}\limits_{
-2A_{p'}(l;n)+1\leq k\leq-A_{p'}(l;n)-1,\atop k \notin\gamma_{2n} }
\Big(1-\frac{|k|}{2A_{p'}(l;n)}\Big)\psi^{p'}(|k|)|k|^{p'-2}\bigg)>
$$
$$
>\pi\bigg(\mathop{\sum}\limits_{
n\leq k\leq A_{p'}(l;n),\atop k \notin\gamma_{2n} }
\psi^{p'}(k)k^{p'-2}+
\mathop{\sum}\limits_{
-A_{p'}(l;n)\leq k\leq-n,\atop k \notin\gamma_{2n} }
\psi^{p'}(|k|)|k|^{p'-2}\bigg)=
$$
\begin{equation}\label{t1}
=
\pi\mathop{\sum}\limits_{
n\leq |k|\leq A_{p'}(l;n),\atop k \notin\gamma_{2n} }
\psi^{p'}(k)k^{p'-2}.
\end{equation}

В силу (\ref{eq3}), (\ref{n1}) i (\ref{t1})
\begin{equation}\label{n2}
I_{1}>\frac{\pi\lambda}{2\Big(\sum\limits_{k=n}^{\infty}\psi^{p'}(k)k^{p'-2}\Big)^{\frac{1}{p}}}\inf\limits_{\gamma_{2n}}\mathop{\sum}\limits_{
n\leq |k|\leq A_{p'}(l;n),\atop k \notin\gamma_{2n} }
\psi^{p'}(k)k^{p'-2}.
\end{equation}

Оскільки при $g_{p}\in\mathfrak{M}_{0}$ функція $\psi^{p'}(t)t^{p'-2}$ монотонно спадає, то
$$
 \inf\limits_{\gamma_{2n}}\mathop{\sum}\limits_{
n\leq |k|\leq A_{p'}(l;n),\atop k \notin\gamma_{2n} }
\psi^{p'}(k)k^{p'-2}=\mathop{\sum}\limits_{
2n\leq |k|\leq A_{p'}(l;n),\atop k \notin\gamma_{2n} }
\psi^{p'}(k)k^{p'-2}=
$$
\begin{equation}\label{eq8}
 =
2\sum\limits_{k=2n}^{A_{p'}(l;n)}
\psi^{p'}(k)k^{p'-2}.
\end{equation}

Покажемо, що за умови, коли функція $g_{s'}(t)=\psi(t)t^{\frac{1}{s'}}$, $1<s<\infty$, ${\frac{1}{s}+\frac{1}{s'}=1}$, така, що $g_{s'}\in\mathfrak{M}_{0}$, то для довільних $l,n\in\mathbb{N}$
\begin{equation}\label{ex4}
\sum\limits_{k=2n}^{A_{s}(l;n)}
\psi^{s}(k)k^{s-2}>\Big(1-\frac{1}{2l}-\frac{s}{\underline{\alpha}_{n}(g_{s'})}\Big)\sum\limits_{k=n}^{\infty}\psi^{s}(k)k^{s-2}.
\end{equation}
Представимо $\sum\limits_{k=2n}^{A_{s}(l;n)}
\psi^{s}(k)k^{s-2}$ у вигляді
\begin{equation}\label{eq88}
\sum\limits_{k=2n}^{A_{s}(l;n)}
\psi^{s}(k)k^{s-2}=\sum\limits_{k=n}^{\infty}
\psi^{s}(k)k^{s-2}-
\sum\limits_{k=n}^{2n-1}\psi^{s}(k)k^{s-2}-
\sum\limits_{k=A_{s}(l;n)+1}^{\infty}\psi^{s}(k)k^{s-2}.
\end{equation}
З (\ref{for3331}) та  спадання функції $\Phi_{s}(\cdot)$ вигляду (\ref{Phi_s}) випливає оцінка
$$
  \sum\limits_{k=A_{s}(l;n)+1}^{\infty}\psi^{s}(k)k^{s-2}\leq\int\limits_{A_{s}(l;n)}^{\infty}\psi^{s}(t)t^{s-2}dt=
$$
\begin{equation}\label{eq10}
=
\Phi_{s}(A_{s}(l;n))<\frac{1}{2l}\Phi_{s}(n)\leq\frac{1}{2l}\sum\limits_{k=n}^{\infty}\psi^{s}(k)k^{s-2}.
\end{equation}

Знайдемо оцінку зверху для суми $\sum\limits_{k=n}^{2n-1}\psi^{s}(k)k^{s-2}$. Для цього скористаємось лемою 3 з \cite{Serdyuk_Stepaniuk2015}.

{\bf Лема 1.} {\it Нехай $\sum\limits_{k=1}^{\infty}\psi^{s}(k)k^{s-2}<\infty$, $1<s<\infty$, $n\in \mathbb{N}$.
  Тоді, якщо функція $g_{s'}:=\psi(t)t^{\frac{1}{s'}}\in\mathfrak{M}_{0}$, $\frac{1}{p}+\frac{1}{p'}=1$, така, що $g_{s'}\in\mathfrak{M}_{0}$, то виконується нерівність
 \begin{equation}\label{lemma1}
\psi^{s}(n)n^{s-1}\leq\frac{s}{\underline{\alpha}_{n}(g_{s'})}\sum\limits_{k=n}^{\infty}\psi^{s}(k)k^{s-2},
\end{equation}
якщо ж  $g_{s'}\in\mathfrak{M}_{C}$, то має місце співвідношення}
\begin{equation}\label{lemma11_1}
\frac{s}{\overline{\alpha}_{n}(g_{s'})}\cdot\frac{n\underline{\alpha}_{n}(g_{s'})}{s+n\underline{\alpha}_{n}(g_{s'})}
\sum\limits_{k=n}^{\infty}\psi^{s}(k)k^{s-2}
\leq\psi^{s}(n)n^{s-1}\leq\frac{s}{\underline{\alpha}_{n}(g_{s'})}\sum\limits_{k=n}^{\infty}\psi^{s}(k)k^{s-2}.
\end{equation}

Враховуючи, що при $g_{s'}\in\mathfrak{M}_{0}$   функція $\psi^{s}(t)t^{s-2}$ спадає тa використовуючи  нерівність (\ref{lemma1}) леми 1, одержимо
\begin{equation}\label{qq1}
  \sum\limits_{k=n}^{2n-1}
\psi^{s}(k)k^{s-2}\leq \psi^{s}(n)n^{s-1}\leq\frac{s}{\underline{\alpha}_{n}(g_{s'})}\sum\limits_{k=n}^{\infty}\psi^{s}(k)k^{s-2}.
\end{equation}
Із (\ref{eq88}), (\ref{eq10}) i (\ref{qq1}) одержимо нерівність (\ref{ex4}).

Застосовуючи нерівність (\ref{ex4}) при $s=p'$, в силу
 формул (\ref{eq9}), (\ref{n2}) i (\ref{eq8}), для довільних $l\in\mathbb{N}$ отримуємо оцінку
$$
e^{\bot}_{2n}(f^{*})_{\infty}\geq
\frac{\lambda}{3}\Big(1-\frac{1}{2l}-\frac{p'}{\underline{\alpha}_{n}(g_{p})}\Big)\Big(\sum\limits_{k=n}^{\infty}\psi^{p'}(k)k^{p'-2}\Big)^{\frac{1}{p'}}=
$$
\begin{equation}\label{qq33}
 =\frac{1}{3\xi(p)}\Big(\frac{\underline{\alpha}_{n}(g_{p})}{p'+\underline{\alpha}_{n}(g_{p})}\Big)^{\frac{1}{p}}
\Big(1-\frac{1}{2l}-\frac{p'}{\underline{\alpha}_{n}(g_{p})}\Big)\Big(\sum\limits_{k=n}^{\infty}\psi^{p'}(k)k^{p'-2}\Big)^{\frac{1}{p'}}.
\end{equation}

Перейшовши до границі в нерівності (\ref{qq33}) при $l\rightarrow\infty$, отримуємо (\ref{n10}). Із (\ref{w1}) i (\ref{n10}) випливає (\ref{theorem_1}).
Теорему 1 доведено.

Неважко переконатись, що  умовам теореми 1
 задовольняють, наприклад, функції
  \begin{equation}\label{ex11}
 \diamond \ \ \psi(t)=t^{-r}, {\frac{1}{p}<r<1}; \ \ \ \ \ \ \ \ \  \ \ \ \ \ \ \ \ \ \ \ \ \ \ \ \ \ \ \ \ \ \ \ \ \  \ \ \ \ \ \ \ \ \ \ \ \ \ \ \ \
 \ \ \ \ \ \ \ \ \ \ \ \ \ \ \ \ \ \ \ \ \ \ \ \ \ \ \ \ \
  \end{equation}
\begin{equation}\label{func}
  \diamond \ \  \psi(t)={t^{-\frac{1}{p}}\ln^{-\gamma}(t+K)}, \ {\gamma>\frac{1}{p'}},  \ K\geq e^{\gamma p'}-1;  \ \  \ \ \ \ \ \ \ \ \ \ \ \ \ \ \ \ \ \ \ \ \ \ \ \ \ \ \ \ \ \ \ \ \ \ \ \ \ \ \ \ \ \ \ \ \
 \end{equation}
 \begin{equation}\label{ex33}
 \diamond \ \ \psi(t)=t^{-\frac{1}{p}}\ln^{-\gamma} (t\!+\!K_{1})
(\ln\ln (t\!+\!K_{2}))^{-\delta}, \ \gamma\geq\frac{1}{p'},  \delta>\frac{1}{p'}, \ K_{2}\geq K_{1} e^{\max\{(\gamma+\delta)p',e\}}-1.
 \end{equation}

{\bf Теорема 2.} {\it
Нехай $1<s<\infty, \ \sum\limits_{k=1}^{\infty}\psi^{s}(k)k^{s-2}<\infty$, $\frac{1}{s}+\frac{1}{s'}=1$
 a функція ${g_{p}(t)=\psi(t)t^{\frac{1}{p}}}$ така, що
 $$
 g_{s'}\in\mathfrak{M}_{0}
 $$
   i
   $$
   \underline{ \alpha}_{1}(g_{s'})=\inf\limits_{t\geq1}\alpha(g_{s'};t)>s.
   $$
  Тоді для довільних $n\in \mathbb{N}$ i $\beta\in\mathbb{R}$ мають місце співвідношення
 \begin{equation}\label{theorem_2}
 \frac{4}{3}K_{\psi,s'}^{(1)}\Big(\sum\limits_{k=n}^{\infty}\psi^{s}(k)k^{s-2}\Big)^{\frac{1}{s}}
  \leq e^{\bot}_{2n}(L_{\beta,1}^{\psi})_{s} \leq e^{\bot}_{2n-1}(L_{\beta,1}^{\psi})_{s} \leq
K_{\psi,s'}^{(2)}\Big(\sum\limits_{k=n}^{\infty}\psi^{s}(k)k^{s-2}\Big)^{\frac{1}{s}},
\end{equation}
де $K^{(1)}_{\psi,s'}$ і $K^{(2)}_{\psi,s'}$  означаються формулами (\ref{k1}) i (\ref{k2}) відповідно. }

  {\bf Доведення теореми 2.} Згідно з теоремою 1 роботи  \cite[с. 245]{Stepaniuk2014} при виконанні умов $g_{s'}\in\mathfrak{M}_{0}$ i
$\sum\limits_{k=1}^{\infty}\psi^{s}(k)k^{s-2}<\infty$, $1<s<\infty$, $\frac{1}{s}+\frac{1}{s'}=1$,
   $n\in \mathbb{N}$, $\beta\in\mathbb{R}$, має місце оцінка
\begin{equation}\label{d1}
 {\cal E}_{n}(L^{\psi}_{\beta,1})_{s}
\leq
K_{\psi,s'}^{(2)}\Big(\sum\limits_{k=n}^{\infty}\psi^{s}(k)k^{s-2}\Big)^{\frac{1}{s}}.
\end{equation}

Тому, враховуючи (\ref{ineq_comp}) i (\ref{d1}), отримуємо оцінку
\begin{equation}\label{d101}
e^{\bot}_{2n}(L_{\beta,1}^{\psi})_{s} \leq e^{\bot}_{2n-1}(L_{\beta,1}^{\psi})_{s} \leq
K_{\psi,s'}^{(2)}\Big(\sum\limits_{k=n}^{\infty}\psi^{s}(k)k^{s-2}\Big)^{\frac{1}{s}}.
\end{equation}

Залишається показати, що
\begin{equation}\label{t2}
e^{\bot}_{2n}(L_{\beta,1}^{\psi})_{s} \geq \frac{4}{3}
K_{\psi,s'}^{(1)}\Big(\sum\limits_{k=n}^{\infty}\psi^{s}(k)k^{s-2}\Big)^{\frac{1}{s}}.
\end{equation}

При довільному $m\in\mathbb{N}$ покладемо
$$
f_{m}(t)=f_{m}(\psi;\beta;t):=
$$
$$
:=\frac{1}{4\pi}\Big(\sum\limits_{k=1}^{m}\!\psi(k)\cos \Big(kt-\frac{\beta\pi}{2}\Big)+2\!\!\sum\limits_{k=m+1}^{2m-1}\!\!\!\Big(1-\frac{k}{2m}\Big)\psi(k)\cos
\Big(kt-\frac{\beta\pi}{2}\Big)\!\!\Big)=
$$
$$
=\frac{1}{8\pi}\bigg(e^{-i\frac{\beta\pi}{2}}\mathop{\sum}\limits_{
1\leq k\leq m}\psi(k)e^{ikt}+
e^{i\frac{\beta\pi}{2}}\mathop{\sum}\limits_{-m\leq k\leq-1 }\psi(|k|)e^{ikt}+
$$
\begin{equation}\label{eq15}
 +2e^{-i\frac{\beta\pi}{2}}\mathop{\sum}\limits_{
m+1\leq k\leq 2m-1}
\Big(1-\frac{k}{2m}\Big)\psi(k)e^{ikt}+
2e^{i\frac{\beta\pi}{2}}\mathop{\sum}\limits_{
-2m+1\leq k\leq m-1}
\Big(1-\frac{|k|}{2m}\Big)\psi(|k|)e^{ikt}\bigg) .
\end{equation}
В \cite[с. 246--247]{Stepaniuk2014} було встановлено, що $f_{m}\in L^{\psi}_{\beta,1}$ при будь--яких $m\in\mathbb{N}$.

Покажемо, що при $m=A_{s}(l;n)$, де $A_{s}(l;n)$ означається рівністю (\ref{for3331}), має місце нерівність
\begin{equation}\label{d2}
e^{\bot}_{2n}(f_{A_{s}(l;n)})_{s}  \geq\frac{1}{4\xi(s')}\Big(\frac{\underline{\alpha}_{n}(g_{s'})}{s+\underline{\alpha}_{n}(g_{s'})}\Big)^{\frac{1}{s'}}
\Big(1-\frac{1}{2l}-\frac{s}{\underline{\alpha}_{n}(g_{s'})}\Big)\Big(\sum\limits_{k=n}^{\infty}\psi^{s}(k)k^{s-2}\Big)^{\frac{1}{s}}, \ l,n\in\mathbb{N}.
\end{equation}

Покладемо
\begin{equation}\label{eq133}
  I_{2}:=\inf\limits_{\gamma_{2n}}\bigg|\int\limits_{-\pi}^{\pi}(f_{A_{s}(l;n)}(t)-S_{\gamma_{2n}}(f_{A_{s}(l;n)};t))\sum\limits_{k=n}^{\infty}
\psi^{s-1}(k)k^{s-2}\cos\Big(kt-\frac{\beta\pi}{2}\Big)dt\bigg|.
\end{equation}
\noindent Використавши  твердження 3.8.1  роботи \cite[с. 137]{Stepanets1}, запишемо
$$
I_{2}\leq\inf\limits_{\gamma_{2n}}\|f_{A_{s}(l;n)}(t)-S_{\gamma_{2n}}(f_{A_{s}(l;n)};t)\|_{s}\Big\|\sum\limits_{k=n}^{\infty}
\psi^{s-1}(k)k^{s-2}\cos\Big(kt-\frac{\beta\pi}{2}\Big)\Big\|_{s'}=
$$
\begin{equation}\label{q4}
=
e^{\bot}_{2n}(f_{A_{s}(l;n)})_{s}\Big\|\sum\limits_{k=n}^{\infty}
\psi^{s-1}(k)k^{s-2}\cos\Big(kt-\frac{\beta\pi}{2}\Big)\Big\|_{s'}, \ \ 1<s<\infty.
\end{equation}

Згідно з формулою (25) роботи \cite[с. 249]{Stepaniuk2014}
$$
\Big\|\sum\limits_{k=n}^{\infty}
\psi^{s-1}(k)k^{s-2}\cos\Big(kt-\frac{\beta\pi}{2}\Big)\Big\|_{s'}\leq
$$
\begin{equation}\label{q5}
\leq
\xi(s')\Big(1+ \frac{s}{\underline{\alpha}_{n}(g_{s'})}\Big)^{\frac{1}{s'}}
\bigg(\sum\limits_{k=n}^{\infty}
\psi^{s}(k)k^{s-2}\bigg)^{\frac{1}{s'}}, \ \ 1<s<\infty.
\end{equation}

В силу (\ref{eq15}) має місце рівність
$$
f_{A_{s}(l;n)}(t)-S_{\gamma_{2n}}(f_{A_{s}(l;n)};t)=
$$
$$
=\frac{1}{8\pi}\bigg(e^{-i\frac{\beta\pi}{2}}\mathop{\sum}\limits_{
1\leq k\leq A_{s}(l;n),\atop k \notin\gamma_{2n} }\psi(k)e^{ikt}+
e^{i\frac{\beta\pi}{2}}\mathop{\sum}\limits_{-A_{s}(l;n)\leq k\leq-1,\atop k \notin\gamma_{2n} }\psi(|k|)e^{ikt}+
$$
$$
+2e^{-i\frac{\beta\pi}{2}}\mathop{\sum}\limits_{
A_{s}(l;n)+1\leq k\leq 2A_{s}(l;n)-1,\atop k \notin\gamma_{2n}}
\Big(1-\frac{k}{2A_{s}(l;n)}\Big)\psi(k)e^{ikt}+
$$
\begin{equation}\label{n11}
+
2e^{i\frac{\beta\pi}{2}}\mathop{\sum}\limits_{
-2A_{s}(l;n)+1\leq k\leq A_{s}(l;n)-1,\atop k \notin\gamma_{2n}}
\Big(1-\frac{|k|}{2A_{s}(l;n)}\Big)\psi(|k|)e^{ikt}\bigg).
\end{equation}
Крім того,
$$
\sum\limits_{k=n}^{\infty}
\psi^{s-1}(k)k^{s-2}\cos\Big(kt-\frac{\beta\pi}{2}\Big)=
$$
\begin{equation}\label{n3}
=\frac{1}{2}\bigg(e^{-i\frac{\beta\pi}{2}}\sum\limits_{k\geq n}
\psi^{s-1}(k)k^{s-2}e^{ikt}+e^{i\frac{\beta\pi}{2}}\sum\limits_{k\leq -n}
\psi^{s-1}(|k|)|k|^{s-2}e^{ikt}\bigg).
\end{equation}
Використавшии (\ref{int_riv}), (\ref{n11}) i (\ref{n3}),
 одержуємо
 $$
 \int\limits_{-\pi}^{\pi}(f_{A_{s}(l;n)}(t)-S_{\gamma_{2n}}(f_{A_{s}(l;n)};t))\sum\limits_{k=n}^{\infty}
\psi^{s-1}(k)k^{s-2}\cos\Big(kt-\frac{\beta\pi}{2}\Big)dt=
 $$
$$
=\frac{1}{8}
\bigg(\mathop{\sum}\limits_{
n\leq k\leq A_{s}(l;n),\atop k \notin\gamma_{2n} }\psi^{s}(k)k^{s-2}+\mathop{\sum}\limits_{
-A_{s}(l;n)\leq k\leq-n,\atop k \notin\gamma_{2n} }\psi^{s}(|k|)|k|^{s-2}+
$$
$$
+2\mathop{\sum}\limits_{
A_{s}(l;n)+1\leq k\leq 2A_{s}(l;n)-1,\atop k \notin\gamma_{2n}}
\Big(1-\frac{k}{2A_{s}(l;n)}\Big)\psi^{s}(k)k^{s-2}+
$$
$$
+2\mathop{\sum}\limits_{
-2A_{s}(l;n)+1\leq k\leq -A_{s}(l;n)-1,\atop k \notin\gamma_{2n}}
\Big(1-\frac{|k|}{2A_{s}(l;n)}\Big)\psi^{s}(|k|)|k|^{s-2}\bigg)>
$$
$$
>\frac{1}{8}
\bigg(\mathop{\sum}\limits_{
n\leq k\leq A_{s}(l;n),\atop k \notin\gamma_{2n} }\psi^{s}(k)k^{s-2}+\mathop{\sum}\limits_{
-A_{s}(l;n)\leq k\leq-n,\atop k \notin\gamma_{2n} }\psi^{s}(|k|)|k|^{s-2}\bigg)=
$$
\begin{equation}\label{eqq32}
=\frac{1}{8}
\mathop{\sum}\limits_{
n\leq |k|\leq A_{s}(l;n),\atop k \notin\gamma_{2n} }\psi^{s}(k)k^{s-2}.
\end{equation}

Отже, в силу (\ref{eq133}) і (\ref{eqq32})
\begin{equation}\label{n4}
I_{2}>\frac{1}{8}\inf\limits_{\gamma_{2n}}
\mathop{\sum}\limits_{
n\leq |k|\leq A_{s}(l;n),\atop k \notin\gamma_{2n} }\psi^{s}(k)k^{s-2}.
\end{equation}

Враховуючи, що при $g_{s'}\in\mathfrak{M}_{0}$   функція $\psi^{s}(t)t^{s-2}$ спадає, то
\begin{equation}\label{eq32}
\inf\limits_{\gamma_{2n}}
\mathop{\sum}\limits_{
n\leq |k|\leq A_{s}(l;n),\atop k \notin\gamma_{2n} }\psi^{s}(k)k^{s-2}=2{\sum\limits_{k=2n}^{A_{s}(l;n)}}
\psi^{s}(k)k^{s-2}.
\end{equation}

З (\ref{ex4}), (\ref{n4}) i (\ref{eq32}) випливає нерівність
\begin{equation}\label{qq5}
  I_{2}>\frac{1}{4}\Big(1-\frac{1}{2l}-\frac{s}{\underline{\alpha}_{n}(g_{s'})}\Big)\sum\limits_{k=n}^{\infty}\psi^{s}(k)k^{s-2}.
\end{equation}

На підставі формул (\ref{q4}), (\ref{q5}) i (\ref{qq5}) отримуємо (\ref{d2}).

З того, що $f_{A_{s}(l;n)}\in L^{\psi}_{\beta,1}$ випливає
$$
e^{\bot}_{2n}(f_{A_{s}(l;n)})_{s}\geq e^{\bot}_{2n}(L^{\psi}_{\beta,1})_{s}, \ \ l\in\mathbb{N}.
$$
Тоді при  $l\in\infty$ з останньої нерівності і нерівності  (\ref{d2})  отримуємо (\ref{t2}).
Теорему 2 доведено.

Оскільки, згідно зі співвідношенням (\ref{lemma11_1}), у випадку, коли $g_{p}\in\mathfrak{M}_{C}$
$$\sum\limits_{k=n}^{\infty}\psi^{p'}(k)k^{p'-2}\asymp\psi^{p'}(n)n^{p'-1},$$
то з теорем 1 і 2  випливає наступне твердження.

{\bf Наслідок 1. }{\it   Нехай $\sum\limits_{k=1}^{\infty}\psi^{p'}(k)k^{p'-2}<\infty$,
 i
$$
 \underline{ \alpha}_{1}(g_{p})=\inf\limits_{t\geq1}\alpha(g_{p};t)>p',
$$
де $g_{p}(t)=\psi(t)(t)t^{\frac{1}{p}}$, $1<p<\infty$, $\frac{1}{p}+\frac{1}{p'}=1$.
  Тоді, якщо  $g_{p}\in\mathfrak{M}_{0}$, то для довільного $\beta\in\mathbb{R}$
  \begin{equation}\label{s1}
e^{\bot}_{n}(L_{\beta,p}^{\psi})_{\infty}\asymp e^{\bot}_{n}(L_{\beta,1}^{\psi})_{p'}\asymp
\Big(\sum\limits_{k=n}^{\infty}\psi^{p'}(k)k^{p'-2}\Big)^{\frac{1}{p'}},
\end{equation}
якщо ж $g_{p}\in\mathfrak{M}_{C}$, то для довільного $\beta\in\mathbb{R}$ }
\begin{equation}\label{j7}
e^{\bot}_{n}(L_{\beta,p}^{\psi})_{\infty}\asymp e^{\bot}_{n}(L_{\beta,1}^{\psi})_{p'}\asymp\psi(n)n^{\frac{1}{p}}.
\end{equation}

Зауважимо, що  коли $g_{p}\in\mathfrak{M}_{0}$ і
\begin{equation}\label{lim}
\lim\limits_{t\rightarrow\infty}\alpha(g_{p};t)=\infty,
\end{equation}
 то порядкові рівності (\ref{j7})  місця не мають, оскільки в цьому випадку виконується оцінка
$$
\psi(n)n^{\frac{1}{p}}=o\bigg(\Big(\sum\limits_{k=n}^{\infty}\psi^{p'}(k)k^{p'-2}\Big)^{\frac{1}{p'}}\bigg), \ n\rightarrow\infty,
$$
яка є наслідком нерівності (\ref{lemma1}).
 Прикладом функцій $\psi$, які задовольняють умови наслідку 1 і для яких виконується умова (\ref{lim}),  є функції виду (\ref{func}) i (\ref{ex33}).

Застосувавши наслідок 1 до функцій $\psi$ виду (\ref{func}) i (\ref{ex33}), отримаємо наступне твердження.

{\bf Наслідок 2. }{\it Нехай $\psi(t)={t^{-\frac{1}{p}}\ln^{-\gamma}(t+K)}$, ${\gamma>\frac{1}{p'}}$, $K\geq e^{\gamma p'}-1$, $1<p<\infty$, $\frac{1}{p}+\frac{1}{p'}=1$ і  $\beta\in\mathbb{R}$. Тоді}
$$
e^{\bot}_{n}(L_{\beta,p}^{\psi})_{\infty}\asymp e^{\bot}_{n}(L_{\beta,1}^{\psi})_{p'}\asymp\psi(n)n^{\frac{1}{p}}\ln^{\frac{1}{p'}} n, \ \  \ n\in\mathbb{N}\setminus\{1\}.
$$

{\bf Наслідок 3. }{\it Нехай $\psi(t)=t^{-\frac{1}{p}}\ln^{-\frac{1}{p'}} (t\!+\!K_{1})
(\ln\ln (t\!+\!K_{2}))^{-\delta}$,   $\delta>\frac{1}{p'}$,  ${K_{2}\geq K_{1}\geq e^{\max\{(\gamma+\delta)p',e\}}-1}$, $1<p<\infty$, $\frac{1}{p}+\frac{1}{p'}=1$,  $\beta\in\mathbb{R}$ i $n\in\mathbb{N}$. Тоді}
$$
e^{\bot}_{n}(L_{\beta,p}^{\psi})_{\infty}\asymp e^{\bot}_{n}(L_{\beta,1}^{\psi})_{p'}\asymp\psi(n)n^{\frac{1}{p}}(\ln n)^{\frac{1}{p'}}
(\ln\ln n)^{\frac{1}{p'}}, \ \  \ n\in\mathbb{N}\setminus\{1,2\}.
$$

{\bf Теорема 3.} {\it  Нехай $\sum\limits_{k=1}^{\infty}\psi(k)<\infty$, а функція  $g(t)=\psi(t)t$ така, що
$$g\in\mathfrak{M}_{0}$$
 i
$$
\underline{\alpha}_{1}(g)=\inf\limits_{t\geq1}\alpha(g;t)>1.
$$

   Тоді, якщо $\cos\frac{\beta\pi}{2}\neq0$, $\beta\in\mathbb{R}$, то для довільного $n\in \mathbb{N}$}
 \begin{equation}\label{theorem_3}
 \frac{1}{12\pi}\Big|\cos\frac{\beta\pi}{2}\Big|\Big(1-\frac{1}{\underline{\alpha}_{1}(g)}\Big)
\sum\limits_{k=n}^{\infty}\psi(k)\leq
e^{\bot}_{2n}(L_{\beta,1}^{\psi})_{\infty}\leq
e^{\bot}_{2n-1}(L_{\beta,1}^{\psi})_{\infty}\leq\frac{1}{\pi}\sum\limits_{k=n}^{\infty}\psi(k).
\end{equation}

 {\bf Доведення теореми 3.}
 В силу теореми 2 роботи \cite[с. 255]{Stepaniuk2014} за умови ${\sum\limits_{k=1}^{\infty}\psi(k)<\infty}$ справедлива нерівність
\begin{equation}\label{j3}
 {\cal E}_{n}(L^{\psi}_{\beta,1})_{\infty}\leq\frac{1}{\pi}\sum\limits_{k=n}^{\infty}\psi(k).
\end{equation}
Із (\ref{ineq_comp}) i (\ref{j3}) маємо
\begin{equation}\label{t4}
 e^{\bot}_{2n}(L_{\beta,1}^{\psi})_{\infty}\leq
e^{\bot}_{2n-1}(L_{\beta,1}^{\psi})_{\infty}\leq\frac{1}{\pi}\sum\limits_{k=n}^{\infty}\psi(k).
\end{equation}

Знайдемо  оцінку знизу величини $e^{\bot}_{2n}(L_{\beta,1}^{\psi})_{\infty}$.

 Покладемо
 $$
 \Psi(x):=\int\limits_{x}^{\infty}\psi(t)dt,
 $$
 \begin{equation}\label{for31}
D(l;n)=D(\psi;l;n):=\big[\Psi^{-1}\big(\frac{1}{2l}\Psi(n)\big)\big]+2n, \ \ l,n\in\mathbb{N},
\end{equation}
i
\begin{equation}\label{eq13}
  I_{3}:=\inf\limits_{\gamma_{2n}}\bigg|\int\limits_{-\pi}^{\pi}(f_{D(l;n)}(t)-S_{\gamma_{2n}}(f_{D(l;n)};t))V_{D(l;n)}(t)dt\bigg|,
\end{equation}
де функція $f_{D(l;n)}(t)$ означається формулою (\ref{eq15}) при $m=D(l;n)$.

Використовуючи
  твердження Д.1.1 з  \cite[с. 391]{Korn}  та формулу (\ref{eq61}), можемо записати оцінку
  $$
  I_{3}\leq\inf\limits_{\gamma_{2n}}\|f_{D(l;n)}(t)-S_{\gamma_{2n}}(f_{D(l;n)};t)\|_{\infty}\|V_{D(l;n)}\|_{1}=
  $$
  \begin{equation}\label{n12}
 =e^{\bot}_{2n}(f_{D(l;n)})_{\infty}\|V_{D(l;n)}\|_{1}\leq3\pi e^{\bot}_{2n}(f_{D(l;n)})_{\infty}.
\end{equation}

Згідно з (\ref{eq15})
$$
f_{D(l;n)}(t)-S_{\gamma_{2n}}(f_{D(l;n)};t)=
$$
$$
=\frac{1}{8\pi}\bigg(e^{-i\frac{\beta\pi}{2}}\mathop{\sum}\limits_{
1\leq k\leq D(l;n),\atop k \notin\gamma_{2n} }\psi(k)e^{ikt}+e^{i\frac{\beta\pi}{2}}\mathop{\sum}\limits_{
-D(l;n)\leq k\leq-1,\atop k \notin\gamma_{2n} }\psi(|k|)e^{ikt}+
$$
$$
+2e^{-i\frac{\beta\pi}{2}}\mathop{\sum}\limits_{
D(l;n)+1\leq k\leq 2D(l;n)-1,\atop k \notin\gamma_{2n}}
\Big(1-\frac{k}{2D(l;n)}\Big)\psi(k)e^{ikt}+
$$
\begin{equation}\label{eq73}
+2e^{i\frac{\beta\pi}{2}}\mathop{\sum}\limits_{
-2D(l;n)+1\leq k\leq -D(l;n)-1,\atop k \notin\gamma_{2n}}
\Big(1-\frac{|k|}{2D(l;n)}\Big)\psi(|k|)e^{ikt}\bigg).
\end{equation}

Із (\ref{val_pus_complex}) при $m=D(l;n)$ маємо
$$
V_{D(l;n)}(t)=\frac{1}{2}\Big(1+\sum\limits_{1\leq k\leq D(l;n)}e^{ikt}+\sum\limits_{-D(l;n)\leq k\leq-1}e^{ikt}+
$$
$$
+2\sum\limits_{D(l;n)+1\leq k\leq 2D(l;n)-1}\Big(1-\frac{k}{2D(l;n)}\Big)e^{ikt}+
$$
\begin{equation}\label{val_pus_complexD}
+2\sum\limits_{-2D(l;n)+1\leq k\leq -D(l;n)-1}\Big(1-\frac{|k|}{2D(l;n)}\Big)e^{ikt}\Big).
\end{equation}

Із (\ref{int_riv}),  (\ref{eq73}) i (\ref{val_pus_complexD})  випливає
$$
\bigg|\int\limits_{-\pi}^{\pi}(f_{D(l;n)}(t)-S_{\gamma_{2n}}(f_{D(l;n)};t))V_{D(l;n)}(t)dt\bigg|=
$$
$$
= \frac{1}{8} \bigg|e^{-i\frac{\beta\pi}{2}}\mathop{\sum}\limits_{
1\leq k\leq D(l;n),\atop k \notin\gamma_{2n} }\psi(k)+e^{i\frac{\beta\pi}{2}}\mathop{\sum}\limits_{
-D(l;n)\leq k\leq-1,\atop k \notin\gamma_{2n} }\psi(|k|)+
$$
$$
+2e^{-i\frac{\beta\pi}{2}}\mathop{\sum}\limits_{
D(l;n)+1\leq k\leq 2D(l;n)-1,\atop k \notin\gamma_{2n}}
\Big(1-\frac{k}{2D(l;n)}\Big)^{2}\psi(k)+
$$
$$
+2e^{i\frac{\beta\pi}{2}}\mathop{\sum}\limits_{
-2D(l;n)+1\leq k\leq -D(l;n)-1,\atop k \notin\gamma_{2n}}
\Big(1-\frac{|k|}{2D(l;n)}\Big)^{2}\psi(|k|)\bigg|=
$$
$$
=\frac{1}{8}\Big|\cos\frac{\beta\pi}{2}\Big|\bigg(\mathop{\sum}\limits_{
1\leq k\leq D(l;n),\atop k \notin\gamma_{2n} }\psi(k)+\mathop{\sum}\limits_{
-D(l;n)\leq k\leq-1,\atop k \notin\gamma_{2n} }\psi(|k|)+
$$
$$
+2\mathop{\sum}\limits_{
D(l;n)+1\leq k\leq 2D(l;n)-1,\atop k \notin\gamma_{2n}}
\Big(1-\frac{k}{2D(l;n)}\Big)^{2}\psi(k)+
$$
$$
+2\mathop{\sum}\limits_{
-2D(l;n)+1\leq k\leq -D(l;n)-1,\atop k \notin\gamma_{2n}}
\Big(1-\frac{|k|}{2D(l;n)}\Big)^{2}\psi(|k|)\bigg)+
$$
$$
+i\sin\frac{\beta\pi}{2}\bigg(-\mathop{\sum}\limits_{
1\leq k\leq D(l;n),\atop k \notin\gamma_{2n} }\psi(k)+\mathop{\sum}\limits_{
-D(l;n)\leq k\leq-1,\atop k \notin\gamma_{2n} }\psi(|k|)-
$$
$$
-2\mathop{\sum}\limits_{
D(l;n)+1\leq k\leq 2D(l;n)-1,\atop k \notin\gamma_{2n}}
\Big(1-\frac{k}{2D(l;n)}\Big)^{2}\psi(k)+
$$
$$
+2\mathop{\sum}\limits_{
-2D(l;n)+1\leq k\leq -D(l;n)-1,\atop k \notin\gamma_{2n}}
\Big(1-\frac{|k|}{2D(l;n)}\Big)^{2}\psi(|k|)\bigg)\bigg|\geq
$$
$$
\geq\frac{1}{8}\Big|\cos\frac{\beta\pi}{2}\Big|\bigg(\mathop{\sum}\limits_{
1\leq k\leq D(l;n),\atop k \notin\gamma_{2n} }\psi(k)+\mathop{\sum}\limits_{
-D(l;n)\leq k\leq-1,\atop k \notin\gamma_{2n} }\psi(|k|)+
$$
$$
+2\mathop{\sum}\limits_{
D(l;n)+1\leq k\leq 2D(l;n)-1,\atop k \notin\gamma_{2n}}
\Big(1-\frac{k}{2D(l;n)}\Big)^{2}\psi(k)+
$$
$$
+2\mathop{\sum}\limits_{
-2D(l;n)+1\leq k\leq -D(l;n)-1,\atop k \notin\gamma_{2n}}
\Big(1-\frac{|k|}{2D(l;n)}\Big)^{2}\psi(|k|)\bigg)>
$$
\begin{equation}\label{t6}
> \frac{1}{8} \Big|\cos\frac{\beta\pi}{2}\Big|\bigg(\mathop{\sum}\limits_{
1\leq k\leq D(l;n),\atop k \notin\gamma_{2n} }\psi(k)+\mathop{\sum}\limits_{
-D(l;n)\leq k\leq-1,\atop k \notin\gamma_{2n} }\psi(|k|)\bigg).
\end{equation}

На підставі (\ref{eq13}) і (\ref{t6}) отримуємо оцінку
$$
 I_{3}>\frac{1}{8}\Big|\cos\frac{\beta\pi}{2}\Big|\inf\limits_{\gamma_{2n}}
 \bigg(\mathop{\sum}\limits_{
1\leq k\leq D(l;n),\atop k \notin\gamma_{2n} }\psi(k)+\mathop{\sum}\limits_{
-D(l;n)\leq k\leq-1,\atop k \notin\gamma_{2n} }\psi(|k|)\bigg)=
$$
$$
=\frac{1}{8}\Big|\cos\frac{\beta\pi}{2}\Big|\inf\limits_{\gamma_{2n}}
 \mathop{\sum}\limits_{
1\leq |k|\leq D(l;n),\atop k \notin\gamma_{2n} }\psi(k)=
\frac{1}{4}\Big|\cos\frac{\beta\pi}{2}\Big|
\sum\limits_{k=n+1}^{D(l;n)}\psi(k)=
$$
\begin{equation}\label{eq16}
=\frac{1}{4}\Big|\cos\frac{\beta\pi}{2}\Big|
\bigg(\sum\limits_{k=n}^{\infty}\psi(k)-\psi(n)-\sum\limits_{k=D(l;n)+1}^{\infty}\psi(k)\bigg).
\end{equation}
 З (\ref{for31}) випливає, що для довільних $ l\in \mathbb{N}$
 \begin{equation}\label{for401}
 \sum\limits_{k=D(l;n)+1}^{\infty}\psi(k)\leq\int\limits_{D(l;n)}^{\infty}\psi(t)dt=\Psi(D(l;n))<
 \frac{1}{2l}\Psi(n)\leq\frac{1}{2l}\sum\limits_{k=n}^{\infty}\psi(k).
 \end{equation}

Далі нам буде корисним наступне твердження роботи \cite[с.
259]{Stepaniuk2014}.

{\bf Лема 2.} {\it Нехай  $\sum\limits_{k=1}^{\infty}\psi(k)<\infty$.
  Тоді, якщо функція $g(t)=\psi(t)t$ така, що
   $g\in\mathfrak{M}_{0}$, то для довільних $n\in \mathbb{N}$
 \begin{equation}\label{lemma_3}
\psi(n)n\leq\frac{1}{\underline{\alpha}_{n}(g)}\sum\limits_{k=n}^{\infty}\psi(k).
\end{equation}

Якщо ж  $g\in\mathfrak{M}_{C}$, то}
\begin{equation}\label{lemma_3_1}
\frac{1}{\overline{\alpha}_{n}(g)}\cdot\frac{n\underline{\alpha}_{n}(g)}{1+n\underline{\alpha}_{n}(g)}\sum\limits_{k=n}^{\infty}\psi(k)
\leq\psi(n)n\leq\frac{1}{\underline{\alpha}_{n}(g)}\sum\limits_{k=n}^{\infty}\psi(k).
\end{equation}

В силу формул (\ref{eq16})--(\ref{lemma_3}), маємо
\begin{equation}\label{qq7}
   I_{3}>\frac{1}{4}\Big|\cos\frac{\beta\pi}{2}\Big|\Big(1-\frac{1}{\underline{\alpha}_{n}(g)n}-\frac{1}{2l}\Big)\sum\limits_{k=n}^{\infty}\psi(k), \ \ l\in\mathbb{N}.
\end{equation}

З (\ref{n12}) та (\ref{qq7})  за умови  $\cos\frac{\beta\pi}{2}\neq0$ отримаємо
$$
e^{\bot}_{2n}(L_{\beta,1}^{\psi})_{\infty}\geq  e^{\bot}_{2n}(f_{D(l;n)})_{\infty}\geq\frac{1}{3\pi}I_{3}>
$$
\begin{equation}\label{eq24}
>\frac{1}{12\pi}\Big|\cos\frac{\beta\pi}{2}\Big|
\Big(1-\frac{1}{\underline{\alpha}_{1}(g)n}-\frac{1}{2l}\Big)\sum\limits_{k=n}^{\infty}\psi(k).
\end{equation}
Перейшовши в формулі  (\ref{eq24}) до границі при $l\rightarrow\infty$, одержимо
\begin{equation}\label{qq8}
e^{\bot}_{2n}(L^{\psi}_{\beta,1})_{\infty}\geq\frac{1}{12\pi}\Big|\cos\frac{\beta\pi}{2}\Big|\Big(1-\frac{1}{\underline{\alpha}_{1}(g)n}\Big)
\sum\limits_{k=n}^{\infty}\psi(k).
\end{equation}

Oб'єднуючи  (\ref{t4}) і (\ref{qq8}) отримуємо (\ref{theorem_3}). Теорему 3 доведено.

{\bf Теорема 4.} {\it
Нехай $\sum\limits_{k=1}^{\infty}\psi(k)<\infty$, а функція  $g(t)=\psi(t)t$, така, що
$$g\in\mathfrak{M}_{0}$$
 i
 \begin{equation}\label{cond_alpha}
\underline{\alpha}_{1}(g)=\inf\limits_{t\geq1}\alpha(g;t)>1.
\end{equation}
   Тоді, якщо $\cos\frac{\beta\pi}{2}=0, \ \beta\in\mathbb{R}$, то для довільних  $n\in \mathbb{N}$
мають місце нерівності
    \begin{equation}\label{theorem_4}
\frac{1}{60\pi}\Big(1-\frac{1}{\underline{\alpha}_{1}(g)}\Big)\psi(n)n \leq e^{\bot}_{2n}(L^{\psi}_{\beta,1})_{\infty}\leq e^{\bot}_{2n-1}(L^{\psi}_{\beta,1})_{\infty}\leq\Big(1+\frac{2}{\pi }\Big)\psi(n)n.
\end{equation}
 }
  {\bf Доведення теореми 4.}
В силу теореми 4 з \cite[с. 262]{Stepaniuk2014} при  виконанні умов ${\sum\limits_{k=1}^{\infty}\psi(k)<\infty}$,
 $g\in\mathfrak{M}_{0}$, $\cos\frac{\beta\pi}{2}=0, \ \beta\in\mathbb{R}$,
справедлива оцінка
\begin{equation}\label{w3}
  e^{\bot}_{2n}(L^{\psi}_{\beta,1})_{\infty}\leq e^{\bot}_{2n-1}(L^{\psi}_{\beta,1})_{\infty}\leq{\cal E}_{n}(L^{\psi}_{\beta,1})_{\infty}\leq\Big(1+\frac{2}{\pi }\Big)\psi(n)n.
\end{equation}
Оцінимо  знизу величину $e^{\bot}_{2n}(L_{\beta,1}^{\psi})_{\infty}$.
Розглянемо функцію
\begin{equation}\label{eq26}
f^{*}_{n}(t)=f^{*}_{n}(\psi;t):=\frac{1}{5\pi n}\Big(\sum\limits_{k=1}^{n}k\psi(k)\cos kt+
\sum\limits_{k=n+1}^{2n}(2n+1-k)\psi(k)\cos kt\Big).
\end{equation}
В \cite[с. 263--265]{Stepaniuk2014} було показано, що $f^{*}_{n}$ належить класу
 $L^{\psi}_{\beta,1}$. Доведемо, що
 \begin{equation}\label{e}
e^{\bot}_{2n}(f^{*}_{n})_{\infty}\geq\frac{1}{60\pi}\Big(1-\frac{1}{\underline{\alpha}_{1}(g)}\Big)\psi(n)n.
\end{equation}

Покладемо
\begin{equation}\label{eq25}
  I_{4}:=\inf\limits_{\gamma_{2n}}\bigg|\int\limits_{-\pi}^{\pi}(f^{*}_{n}(t)-S_{\gamma_{2n}}(f^{*}_{n};t))V_{2n}(t)dt\bigg|,
\end{equation}
де  $V_{m}$ --- суми Валле Пуссена вигляду (\ref{val_pus2}).

Використавши   твердження Д.1.1 з  \cite[с. 391]{Korn}  та нерівність (\ref{eq61}), отримаємо
\begin{equation}\label{eq27}
I_{4}\leq\inf\limits_{\gamma_{2n}}\|f^{*}_{n}(t)-S_{\gamma_{2n}}(f^{*}_{n};t)\|_{\infty}\|V_{2n}\|_{1}
\leq3\pi \ e^{\bot}_{2n}(f_{n}^{*})_{\infty}.
\end{equation}

Оскільки, в силу формули   (\ref{eq26}) має місце рівність
$$
f^{*}_{n}(t)-S_{\gamma_{2n}}(f^{*}_{n};t)=\frac{1}{10\pi n}\Big(\mathop{\sum}\limits_{
1\leq k\leq n,\atop k \notin\gamma_{2n} }k\psi(k)e^{ikt}+\mathop{\sum}\limits_{
-n\leq k\leq -1,\atop k \notin\gamma_{2n} }|k|\psi(|k|)e^{ikt}+
$$
$$
+\mathop{\sum}\limits_{
n+1\leq k\leq 2n,\atop k \notin\gamma_{2n} }(2n+1-k)\psi(k)e^{ikt}+
\mathop{\sum}\limits_{
-2n\leq k\leq -n-1,\atop k \notin\gamma_{2n} }(2n+1-k)\psi(|k|)e^{ikt}
\Big)
$$
а в силу (\ref{val_pus_complex}) --- рівність
$$
V_{2n}(t)=\frac{1}{2}\Big(1+\sum\limits_{1\leq k\leq 2n}e^{ikt}+\sum\limits_{-2n\leq k\leq -1}e^{ikt}+2\sum\limits_{2n+1\leq k\leq 4n-1}\Big(1-\frac{k}{2n}\Big)e^{ikt}
+
$$
$$
+2\sum\limits_{-4n+1\leq k\leq-2n-1}\Big(1-\frac{|k|}{2n}\Big)e^{ikt}\Big),
$$
то застосовуючи формули  (\ref{int_riv}), знаходимо
$$
\int\limits_{-\pi}^{\pi}(f^{*}_{n}(t)-S_{\gamma_{2n}}(f^{*}_{n};t))V_{2n}(t)dt=
$$
$$
=\frac{1}{10 n}\bigg(
 \mathop{\sum}\limits_{
1\leq k\leq n,\atop k \notin\gamma_{2n} }k\psi(k)+\mathop{\sum}\limits_{
-n\leq k\leq -1,\atop k \notin\gamma_{2n} }|k|\psi(|k|)+
$$
\begin{equation}\label{t7}
+\mathop{\sum}\limits_{
n+1\leq k\leq 2n,\atop k \notin\gamma_{2n} }(2n+1-k)\psi(k)+
\mathop{\sum}\limits_{
-2n\leq k\leq -n-1,\atop k \notin\gamma_{2n}}(2n+1-|k|)\psi(|k|)
\bigg).
\end{equation}
 Враховуючи формули (\ref{eq25}) і (\ref{t7}), монотонне спадання функції $g$, та виконуючи елементарні перетворення, запишемо оцінку величини $I_{4}$
$$
 I_{4}=\frac{1}{10 n}\inf\limits_{\gamma_{2n}}\bigg(
 \mathop{\sum}\limits_{
1\leq k\leq n,\atop k \notin\gamma_{2n} }k\psi(k)+\mathop{\sum}\limits_{
-n\leq k\leq -1,\atop k \notin\gamma_{2n} }|k|\psi(|k|)+
$$
$$
+\mathop{\sum}\limits_{
n+1\leq k\leq 2n,\atop k \notin\gamma_{2n} }(2n+1-k)\psi(k)+
\mathop{\sum}\limits_{
-2n\leq k\leq -n-1,\atop k \notin\gamma_{2n}}(2n+1-|k|)\psi(|k|)
\bigg)>
$$
$$
>\frac{1}{5 n}\sum\limits_{k=n+1}^{2n}\psi(k)(2n+1-k)\geq\frac{\psi(2n)}{5 n}\sum\limits_{k=n+1}^{2n}(2n+1-k)=
$$
\begin{equation}\label{eq28}
=\psi(2n)\frac{n+1}{10}>\frac{1}{10}\psi(2n)n.
\end{equation}
Використавши  співвідношення
(\ref{eq27}) і (\ref{eq28}), отримаємо
\begin{equation}\label{eq288}
e^{\bot}_{2n}(L^{\psi}_{\beta,1})_{\infty}\geq e^{\bot}_{2n}(f^{*}_{n})_{\infty}\geq\frac{1}{3\pi}I_{4}\geq
\frac{1}{30\pi} \psi(n)n\frac{\psi(2n)}{\psi(n)}
=\frac{1}{60\pi}\psi(n)n\frac{g(2n)}{g(n)}.
\end{equation}
Оскільки, як показано в   \cite[с. 266]{Stepaniuk2014} за умови (\ref{cond_alpha}) виконується нерівність
$$
\frac{g(2n)}{g(n)}> 1-\frac{1}{\underline{\alpha}_{1}(g)},
$$
то з
 (\ref{eq288}) випливає оцінка (\ref{e}).

Із (\ref{w3}) і (\ref{e}) випливає  (\ref{theorem_4}).
Теорему 4 доведено.

Оскільки $g\in\mathfrak{M}_{0}$, де $g(t)=\psi(t)t$, то згідно з \cite[с.
175]{Stepanets1} виконується нерівність $\frac{g(2n)}{g(n)}>K_{1}$.
 Тоді з  (\ref{eq288}) отримуємо оцінку
\begin{equation}\label{qq8811}
e^{\bot}_{2n}(L^{\psi}_{\beta,1})_{\infty}\geq K_{2}
\psi(n)n, \ \ \cos\frac{\beta\pi}{2}=0, \ \ \ \beta\in\mathbb{R}.
\end{equation}

Крім того, очевидно, що при досить великих $n$ справджується нерівність ${\underline{\alpha}_{1}(g)n>K_{3}>1}$.
  Тоді з (\ref{qq8}) маємо
\begin{equation}\label{qq88}
e^{\bot}_{2n}(L^{\psi}_{\beta,1})_{\infty}\geq K_{4}
\sum\limits_{k=n}^{\infty}\psi(k), \ \ \cos\frac{\beta\pi}{2}\neq0, \ \ \ \beta\in\mathbb{R}.
\end{equation}

Згідно зі співвідношенням (\ref{lemma_3_1}) леми 2, якщо $g\in\mathfrak{M}_{C}$, то
\begin{equation}\label{qq881}
\sum\limits_{k=n}^{\infty}\psi(k)\asymp\psi(n)n.
\end{equation}
Із (\ref{t4}), (\ref{w3}), (\ref{qq8811})--(\ref{qq881}) приходимо до наступного  твердження.

{\bf Теорема 5. }{\it   Нехай $\sum\limits_{k=1}^{\infty}\psi(k)<\infty$ і  $\beta\in\mathbb{R}$ .
  Тоді, якщо функція $g(t)=\psi(t)t$, така, що
  $$
  g\in\mathfrak{M}_{0},
  $$
   то
\begin{equation}\label{co}
  e^{\bot}_{n}(L^{\psi}_{\beta,1})_{\infty}\asymp {\left\{\begin{array}{cc}
\sum\limits_{k=n}^{\infty}\psi(k), \ \ & \cos\frac{\beta\pi}{2}\neq0, \\
\psi(n)n,  & \cos\frac{\beta\pi}{2}=0, \
  \end{array} \right.}
 \end{equation}
 якщо ж $g\in\mathfrak{M}_{C}$, то  }
\begin{equation}\label{co1}
   e^{\bot}_{n}(L^{\psi}_{\beta,1})_{\infty}\asymp \psi(n)n.
\end{equation}

Неважко переконатись, що умови теореми 5 задовольняють, наприклад,  функції:
\begin{equation}\label{ex1}
\diamond \ \  \psi(t)=t^{-r}, \ \ r>1; \ \ \ \ \ \ \ \ \ \ \ \ \ \ \ \ \ \ \ \ \ \ \ \ \ \ \ \ \ \ \ \ \ \ \ \ \ \
\ \ \ \ \ \ \ \ \ \ \ \ \ \ \ \ \ \ \ \ \ \ \ \ \ \ \ \ \ \ \ \ \ \ \ \ \ \ \ \ \ \ \ \ \ \ \ \ \ \ \
 \end{equation}
\begin{equation}\label{ex2}
 \diamond \ \ \psi(t)={t^{-1}\ln^{-\gamma}(t+K)}, \ K>0, \ {\gamma>1};
 \ \ \ \ \ \ \ \ \ \ \ \ \ \ \ \ \ \ \ \ \ \ \ \ \ \ \ \ \ \ \ \ \ \ \ \ \ \ \ \ \ \ \ \ \ \ \ \ \ \ \ \ \ \ \ \ \ \ \ \ \ \
 \end{equation}
\begin{equation}\label{ex3}
 \diamond \ \ \psi(t)={t^{-1}\ln^{-\gamma} (t+K_{1})
(\ln\ln (t+K_{2}))^{-\delta}}, \gamma\geq1, \delta>1, K_{1}>0,  K_{2}> e-1. \ \ \ \ \ \ \ \ \ \ \ \ \ \ \ \ \ \ \ \
 \end{equation}

Зауважимо, що  коли $g\in\mathfrak{M}_{0}$, $g(t)=\psi(t)t$ і
\begin{equation}\label{lim1}
\lim\limits_{t\rightarrow\infty}\alpha(g;t)=\infty,
\end{equation}
 то  в цьому випадку  виконується оцінка
$$
\psi(n)n=o\bigg(\sum\limits_{k=n}^{\infty}\psi(k)\bigg), \ n\rightarrow\infty,
$$
яка є наслідком нерівності (\ref{lemma_3}).

 Прикладом функцій $\psi(t)$, які задовольняють умови теореми 5 і для яких виконується умова (\ref{lim1}),  є функції виду (\ref{ex2}) та (\ref{ex3}).

Наведемо порядкові оцінки величин $e^{\bot}_{n}(L^{\psi}_{\beta,1})_{\infty}$ для функцій  виду (\ref{ex1})--(\ref{ex3}).

{\bf Наслідок 4. }{\it Нехай  $\psi(t)=t^{-r}$,  $r>1$ і  $\beta\in\mathbb{R}$. Тоді}
$$
 e^{\bot}_{n}(W^{r}_{\beta,1})_{\infty}\asymp n^{-r+1}.
$$

{\bf Наслідок 5. }{\it Нехай  $\psi(t)=t^{-1}\ln^{-\gamma} (t+K)$,  $\gamma>1$, $K>0$,  $\beta\in\mathbb{R}$,
 i $n\in\mathbb{N}\setminus\{1\}$. Тоді}
 $$
 e^{\bot}_{n}(L^{\psi}_{\beta,1})_{\infty}\asymp {\left\{\begin{array}{cc}
\psi(n)n\ln n, \ \ & \cos\frac{\beta\pi}{2}\neq0, \\
\psi(n)n,  \ \ \ \ \  & \cos\frac{\beta\pi}{2}=0. \
  \end{array} \right.}
 $$

{\bf Наслідок 6. }{\it Нехай  $\psi(t)={t^{-1}\ln^{-1} (t+K_{1})
(\ln\ln (t+K_{2}))^{-\delta}}$,  $\delta>1$, ${K_{1}>0},\  {K_{2}>e-1}$,  $\beta\in\mathbb{R}$,
 i $n\in\mathbb{N}\setminus\{1,2\}$. Тоді}
 $$
 e^{\bot}_{n}(L^{\psi}_{\beta,1})_{\infty}\asymp {\left\{\begin{array}{cc}
\psi(n)n\ln n \ln(\ln n), \ \ & \cos\frac{\beta\pi}{2}\neq0, \\
\psi(n)n,  \ \ \ \ \  & \cos\frac{\beta\pi}{2}=0. \
  \end{array} \right.}
 $$

E-mail: \href{mailto:serdyuk@imath.kiev.ua}{serdyuk@imath.kiev.ua},
\href{mailto:tania_stepaniuk@ukr.net}{tania$_{-}$stepaniuk@ukr.net}

\end{document}